\documentclass[a4paper,12pt]{article}
\usepackage[cp1251]{inputenc}
\usepackage{amssymb,amsmath,amsthm}
\usepackage[russian]{babel}
\usepackage{amssymb,amsmath,amsthm}
\usepackage{indentfirst}
\usepackage{hyperref,hypbmsec}
\usepackage[pdftex]{graphicx}
\numberwithin{equation}{section}
\newtheorem{Th}{\hskip\parindent Теорема}[section]
\newtheorem{Le}{\hskip\parindent Лемма}[section]
\newtheorem{Sl}{\hskip\parindent Следствие}[section]
\newtheorem{Zam}{\hskip\parindent Замечание}[section]
\newtheorem{Hyp}{\hskip\parindent Гипотеза}[section]

\newcommand{\A}{\mathcal{A}}
\newcommand{\B}{\mathcal{B}}

\newcommand{\R}{\mathfrak{R}}

\newcommand{\D}{\mathfrak{D}}
\newcommand{\N}{\mathbb{N}}
\newcommand{\M}{\mathfrak{M}}
\newcommand{\NN}{\mathfrak{N}}
\newcommand{\1}{\mathbf{1}}

\newcommand{\q}{\mathbf{q}}

\newcommand{\rr}{\mathbb{R}}
\newcounter{propet}

\renewcommand{\le}{\leqslant}\renewcommand{\ge}{\geqslant}

\begin{document}

\author{
И. Д. Кан,\footnote{Работа выполнена при поддержке РФФИ (грант  12-01-00681-a)} }

\title{
\begin{flushright}
\small{ Посвящается светлой памяти\\ Н. М. Коробова.}
\end{flushright}
\begin{flushleft}
УДК 511.321 + 511.31
\end{flushleft}
Усиление теоремы Бургейна --- Конторовича - III }
\date{}
\maketitle

\begin{abstract}
В настоящей работе доказывается, что в натуральном ряду чисел имеется положительная пропорция знаменателей тех конечных цепных дробей,  все неполные частные которых принадлежат алфавиту $\{1,2,3,4,10\}$. Ранее аналогичная теорема была известна лишь для алфавита $\{1,2,3,4,5\}$, либо для алфавитов большей мощности. 

\noindent

\textbf
{Библиография:} 14 названий.\\
\noindent

\textbf
{Ключевые слова и выражения:\,} цепная дробь, континуант, тригонометрическая сумма, гипотеза Зарембы. \par

\end{abstract}

\setcounter{Zam}0
\section{История вопроса}
Через $\left     [d_1,d_2,\ldots,d_{k}\right]$ обозначают конечную цепную дробь

\begin{equation}
\label{cont.fraction}
[d_1,d_2,\ldots,d_k]=\frac{1}{d_1+{\atop \ddots {+\frac{1}{d_k}}}}
\end{equation}
с натуральными неполными частными 
$\ d_1,d_2,\ldots,d_{k}$ (где $k$ --- натуральное), а через $\R_{\A}$ --- множество рациональных чисел $\frac{b}{d}$ , представимых  конечными цепными дробями  с неполными частными из некоторого конечного алфавита $\A \subseteq \N$ (всюду далее $|\A|\ge2$ ):

  $$
\R_{\A}=\left\{\frac{b}{d}=[d_1,d_2,\ldots,d_{k}]\Bigl|{} \  d_j\in\A \,\mbox{для}\, j=1,\ldots,k\right\}.
$$
Наконец, через $D_{\A}(N)$   для каждого $N\in \N$ обозначают множество знаменателей $d$  чисел $\frac{b}{d}\in\R_{\A}$, таких что $d$ не превосходит $N$: 

$$\D_{\A}(N)=\left\{d\in \N\Bigl|{} \  \exists b:{}\  (b,d)=1, {}\  \frac{b}{d}\in\R_{\A},{}\ d\le N\right\}.$$

\begin{Hyp}\label{h1.1} (гипотеза Зарембы  \cite{Zaremba}, 1971). Существует константа $A$ (скорее всего, $A=~5$), такая что для любого $N\in \N$ для алфавита 
\begin{equation}
\label{AA}
\A= 1,2,...,A 
\end{equation}
имеет место равенство 
$|\D_{\A}(N)|=N$.\end{Hyp}

Подробный обзор результатов, связанных с гипотезой 1.1, можно найти в  работах \cite{BK},\cite{NG}.  Отметим лишь, что пргопагандистом и энтузиастом этой темы задолго до 1971 года был профессор Н. М. Коробов. Он, в частности, доказал \cite{Korobov} , что для простого  $d$ существует натуральное число $b<d$, такое что $\frac{b}{d}\in\R_{\A}$ для алфавита $\A$ вида~(\ref{AA}) при $A\le \log d$.

Пусть  $\delta_{\A}$ --- хаусдорфова размерность множества бесконечных цепных дробей с неполными частными из произвольного конечного алфавита ${\A}$. Бургейн и Конторович в 2011 году доказали, в частности, следующие две теоремы.

\begin{Th}\label{1.1} \cite[стр.13, теорема 1.25]{BK}. Для каждого
алфавита ${\A}$, такого что
\begin{equation}
 \label{c}                                                                                                                                                                                                             
 \delta_{\A}>\frac{307}{312}=0.9839\ldots,                                                                                                                                   
\end{equation}
справедливо неравенство (``положительная пропорция''):
\begin{equation}
 \label{lc}                                                                                                                                                                                                             
 |\D_{\A}(N)|>> N.                                                                                                                                   
\end{equation}
\end{Th}

\begin{Th}\label{t1.2}
 \cite[стр.13, теорема 1.27]{BK}.  Для каждого алвавита ${\A}$, удовлетворяющего условию (\ref{c}), найдется константа $c=c({\A})>0$, такая что
 
 \begin{equation}
 \label{lld}                                                                                                                                                                                                             
N-|\D_{\A}(N)|{} \ {}\ll N^{1-\frac{c}{\log\log N}};                                                                                                                               
 \end{equation}
кроме того, каждое число $d$ из отрезка $[1,N]$, за исключением  не более чем $\ll N^{1-\frac{c}{\log\log N}}$ чисел, появляется во множестве $\D_{\A}(N)$ с кратностью
 \begin{equation}
 \label{cdd}                                                                                                                                                                                                             
 \left   |\left \{b \Bigl|{}  {}\  \ {} \frac{b}{d}\in\R_{\A},{}\  (b,d)=1\right\} \right| {>>_\epsilon}{}\ N^ {2\delta_{\A}-1-\epsilon}.                                                                                                                                   
\end{equation}
\end{Th}
     Результаты Хенсли \cite{Hen1}   дают веские основания предполагать, что неравенству (\ref{c}) удовлетворяет алфавит (\ref{AA}) при  $A=50$ (но не $A=34$ --- ввиду результатов Дженкинсона  \cite{Jenkinson} ).
     
     Автор  настоящей статьи совместно с  Д.~А. Фроленковым (\cite{FK1} --- \cite{FK4}) усилили  первую из упомянутых теорем  Бургейна  ---  Конторовича.  Усиление состояло в  понижении достаточной нижней оценки   $\delta_{\A}$ вида (\ref{c}) при сохранении итогового неравенства (\ref{lc}). Так, в работе \cite{FK4} неравенство (\ref{lc}) было доказано лишь при условии 
     \begin{equation}
 \label{h}                                                                                                                                                                                                             
 \delta_{\A}>\frac{5}{6}=0.8333\ldots,                                                                                                                                   
\end{equation}
соответствующем алфавиту (\ref{AA}) при  $A=5$, но аналогичное обобщение теоремы \ref{t1.2}  выведено не было. 

Недавно (2013 г.), объединяя методы работ \cite{BK}  и \cite{FK4}, С. Хуанг  \cite{Huang}   доказал, что свойства  (\ref{lld}) и (\ref{cdd}) при выполнении неравенства (\ref{h}) также справедливы. Работа \cite{Huang} содержит два примечательных достижения. Одно из них  состоит в упрощении метода Бургейна --- Конторовича в части вывода неравенств (\ref{lld}) и (\ref{cdd}). Другое, неменее важное, --- может быть сформулировано в виде \textbf{тезиса}: 

\textit{какова бы ни была нижняя грань хаусдорфовой размерности $\delta_{\A}$, позволяющая вывести неравенство  (\ref{lc}) в рамках рассматриваемого метода, та же самая оценка $\delta_{\A}$ позволяет получать утверждения (\ref{lld}) и (\ref{cdd}) аналогичным образом, используя аргументы работы  \cite{Huang}. }

\section{Основные результаты работы}

Основной результат настоящей статьи --- следующий. 

\begin{Th}\label{2.1}
 Для произвольного алфавита $ \A$, такого что 
 \begin{equation}
 \label{hy}                                                                                                                                                                                                             
 \delta_{\A}>\frac{4}{5}=0.8,                                                                                                                                   
\end{equation}
имеет место неравенство $\left    |\D_{\A}(N)\right|>> N$.
\end{Th}

Применение методов Хуанга позволяет теперь вывести следующую теорему.

\begin{Th}\label{2.2}
 Для произвольного алфавита $ \A$, удовлетворяющего неравенству  (\ref{hy}), имеют место оценки  (\ref{lld}) и (\ref{cdd}).
\end{Th}

\begin{Zam}\label{z2.1}
  Согласно результатам работы Дженкинсона \cite{Jenkinson} , неравенству (\ref{hy}) удовлетворяют все алфавиты вида
$$\A= \{1,2,3,4,n\},$$
где число n может принимать любое из значений 6,7,8,9,10. \end{Zam}

\section{Благодарности.} Автор благодарит профессора Н. Г. Мощевитина за постановку задачи и неоднократное обсуждение темы статьи. Также автор благодарен И. Д. Шкредову и И. С. Резвяковой за вопросы и комментарии во время докладов автора. Кроме того, автор весьма благодарен Д. А. Фроленкову за многие  полезные советы --- в частности, совет использовать методы работы  \cite{Huang}. 

\section{Обозначения.} Всюду далее $\epsilon_0~\in~\left  (0,\frac{1}{2500}\right)$  --- произвольно малая положительная константа, участвующая в построении ансамбля $\Omega_N$. Знак Виноградова $f(N)\ll g(N)$ для двух произвольных функций $f(N)$ и $g(N)$ обозначает существование константы  $C,$ зависящей только от $\A$ и $\epsilon_0,$ такой что  $|f(N)|\le Cg(N).$  Если при этом $C=C(\epsilon)$ для произвольного $\epsilon>0$, то используется обозначение $f(N)\ll_{\epsilon} g(N,\epsilon)$.  Также используется традиционное обозначение  $e(x)=\exp(2\pi ix).$ Наибольший общий делитель двух целых чисел $a$ и $b$ обозначается через $(a,b).$ Если $P$ --- утверждение, то ${\1_{\{P\}}}=1$, когда $P$ истинно, и ${\1_{\{P\}}}=0,$ когда  $P$ ложно.  Мощность конечного множества  $S$ обозначается через $|S|$. Для действительного числа $\alpha$ через  $[\alpha]$,$\{\alpha\}$ и $\|\alpha\|$ обозначаются, соответственно,  целая часть от  $\alpha$, дробная доля $\alpha$ и расстояние от $\alpha$ ближайшего целого:
$$[\alpha]=\max\left\{z\in\mathbb{Z}|\,z\le \alpha \right\},{}\ {}\   \{\alpha\}=\alpha-[\alpha],{}\ {}\ 
\|\alpha\|=\min\left\{|z-\alpha|\Bigl|\,z\in\mathbb{Z}\right\}. $$
Кроме того, если $g$ \--- матрица, то $\|g\|$ \--- ее норма (определенная ниже в параграфе \ref{5}).

\section{Континуанты и матрицы.} \label{5}Континуантом натуральных чисел $\ d_1,d_2,\ldots,d_{k}$ называется число $\ \left\langle d_1,d_2,\ldots,d_{k}\right\rangle$, равное знаменателю цепной дроби (\ref{cont.fraction}), несократимому с ее числителем (в том числе, континуант пустой последовательности считается равным одному). Если последовательность $\ d_1,d_2,\ldots,d_{k}$ обозначена через $D$, то через $\stackrel{\leftarrow}{D}$ обозначают последовательность $\ d_{k},d_{k-1},\ldots,d_1,$ а через $D^{-}$ и $D_{-}$ --- последовательности $\ d_1,d_2,\ldots,d_{k-1}$ и $\ d_{2},d_{3},\ldots,d_{k}$, соответственно. Хорошо известно (например,  \cite{knut}), что для произвольных конечных последовательностей $D,{} X$ выполнено неравенство
\begin{equation}\label{continuant inequality}
\langle D\rangle\langle X\rangle\le
\langle D,X\rangle\le
2\langle D\rangle\langle X\rangle,
\end{equation}
следующее из равенства
\begin{equation}\label{continuant equality}
\langle D,X\rangle=
\left(1+\left[\stackrel{\leftarrow}{D}\right][X]\right)\langle D\rangle\langle X\rangle,
\end{equation}
и что для матрицы
\begin{gather}\label{matrix-def}
g=
\begin{pmatrix}
a & b \\
c & d
\end{pmatrix}=
\begin{pmatrix}
0 & 1 \\
1 & d_1
\end{pmatrix}
\begin{pmatrix}
0 & 1 \\
1 & d_2
\end{pmatrix}\ldots
\begin{pmatrix}
0 & 1 \\
1 & d_k
\end{pmatrix}
\end{gather}
выполнены соотношения $ (b,d)=(c,d)=1,{}  \ {} a<c<d,{}   \ {}  a<b<d,$
\begin{equation}\label{matrix-no}
\frac{b}{d}=\frac{\ \left\langle D_-\right\rangle}{\left\langle D\right\rangle\ }=[d_1,d_2,\ldots,d_{k}]=[D], {}\\\\\\\\\\\ {}
 \frac{c}{d}=\frac{\ \left\langle D^-\right\rangle}{\ \left\langle D\right\rangle}=[d_{k},d_{k-1},\ldots,d_{1}]=\left[\stackrel{\leftarrow}{D}\right].
 \end{equation}
     
Следовательно,  для нормы такой матрицы $g$ имеет место равенство
\begin{equation}\label{matrix-normm}
\|g\|=\max\{|a|,|b|,|c|,|d|\}=d=\langle d_1,d_2,\ldots,d_{k}\rangle.
\end{equation}
    
  \begin{Le}\label{5.1} (сходное утверждение имеется в   \cite{DKM}). Пусть  $D,T,W$ --- конечные последовательности чисел из алфавита (\ref{AA}), в том числе $T$ и $W$ --- непустые и имеют различие в своих первых элементах и равную длину. Тогда справедливо неравенство
\begin{equation}
\label{norm}
    |[D,T]-[D,W]|\ge\frac{1}{(2A)^{4}\langle D\rangle^2.}
\end{equation}
\end{Le}

Доказательство. Прежде всего, для любой непустой конечной последовательности $V$ в алфавите (\ref{AA}) имеет место оценка 
\begin{equation}
\label{norm0}
    \frac{1}{2A}\le[A,1]\le [V]\le 1.
\end{equation}
  
  Пусть, для начала, последовательность $D$ не пуста. Тогда, если $X$ --- любая из последовательностей $T$ или $W$, то, используя равенства (\ref{continuant equality}) и  (\ref{matrix-no}), получаем:
$$[D,X]=\frac
                      {\langle (D,X)_{-}\rangle}
                      {\langle D,X\rangle}
=\frac
                {\left(1+\left[\left({\stackrel{\leftarrow}{D}}\right)^{-}\right][X]\right)\  \left\langle {D}_{-}\right\rangle \left\langle X\right\rangle}
                {\left(1+\left[\stackrel{\leftarrow}{D}\right][X]\right)\ \left\langle D\right\rangle\ \left\langle X\right\rangle}
=[D]\frac
                {1+[X]\left[\left({\stackrel{\leftarrow}{D}}\right)^{-}\right]}
                {1+[X]\left[{\stackrel{\leftarrow}{D}}\right]}.$$
 Следовательно, обозначая левую часть неравенства  (\ref{norm}) через $\sigma$, получаем:
\begin{equation}
\label{nor}
  \sigma=[D]\frac
                           {\left|\left(\left[{\stackrel{\leftarrow}{D}}\right]-\left[\left({\stackrel{\leftarrow}{D}}\right)^{-}\right]\right)([W]-[T])\right|}
                           {\left(1+[T]\left[{\stackrel{\leftarrow}{D}}\right]\right)\left(1+[W]\left[{\stackrel{\leftarrow}{D}}\right]\right)}
                           \ge\frac{\left|\left(\left[{\stackrel{\leftarrow}{D}}\right]-\left[\left({\stackrel{\leftarrow}{D}}\right)^{-}\right]\right)([W]-[T])\right|}{8A}
\end{equation}
   в виду неравенства (\ref{norm0}), примененного к каждому из случаев  $V=D,\stackrel{\leftarrow}{D},T,W$.  Так как, кроме того, имеют место оценки
     $$\left|\left[{\stackrel{\leftarrow}{D}}\right]
     -\left[\left({\stackrel{\leftarrow}{D}}\right)^{-}\right]\right|=
   \frac{1}{\left\langle D\right\rangle \left\langle {D_{-}}
\right\rangle}\ge\frac{1}{\langle D\rangle^2},$$
\begin{equation}
\label{noryy}
|[W]-[T]|\ge[A-1,1]-[A,A,1]\ge\frac{1}{2A^3},
\end{equation}
то, подставляя их в (\ref{nor}), получаем (\ref{norm}). 

Наконец, если последовательность $D$ пуста, то нужное неравенство  (\ref{norm}) следует непосредственно из  (\ref{noryy}). Лемма доказана.

\section{Основные свойства ансамбля $\Omega_N$.}\label{6} Через $\Gamma_{\A}$ обозначают мультипликативную полугруппу $\Gamma_{\A}\subseteq SL\left(2,\mathbb{Z}\right)$ с единицей  $E=
\begin{pmatrix}
1 & 0 \\
0 & 1
\end{pmatrix}$, порожденную попарными матричными произведениями
\begin{equation*}
\begin{pmatrix}
0 & 1 \\
1 & a
\end{pmatrix}
\begin{pmatrix}
0 & 1 \\
1 & b
\end{pmatrix}=
\begin{pmatrix}
1   & b \\
a & ab+1
\end{pmatrix},
\end{equation*}
где $a,  b \in\A.$ Обозначим через $V_{\A}$ множество слов четной длины (= состоящих из четного количества букв) в алфавите $\A$. Далее всюду будем использовать взаимную однозначность отображения
$\B: V_{\A}\rightarrow\Gamma_{\A},$ определенного формулой

\begin{gather}\label{200norm}
\B(d_1,d_2,\ldots,d_{2k}) =
\begin{pmatrix}
0 & 1 \\
1 & d_1
\end{pmatrix}
\begin{pmatrix}
0 & 1 \\
1 & d_2
\end{pmatrix}\ldots
\begin{pmatrix}
0 & 1 \\
1 & d_{2k}
\end{pmatrix}.
\end{gather}

Напомним, что подмножество $\Theta $ матриц $g\in\Gamma_{\A}$ называется предансамблем  \cite[параграф 8]{FK3}, если для любых двух матриц $g_{1},g_2 \in\Gamma_{\A}$ и любого положительного числа~ $\epsilon$ выполнены соотношения $\|g_{1}\|   \ll  ~  \|g_{2}\| \ll   ~ \|g_{1}\|$ и $| \Theta |>>_{\epsilon} \left\|g_{1}\right\| ^{2\delta_{\A}-\epsilon}$.  По произвольному достаточно большому числу $N$ и по малому параметру $\epsilon_0\in \left(0,\frac{1}{2500}\right)$ в   \cite{FK3} была построена конечная последовательность
\begin{equation}
\label{hgfj}
\left\{N_{-J-1},N_{-J},\ldots,N_{-1},N_0,N_1,\ldots,N_{J+1}\right\},
\end{equation}
где $N_{J+1}=N,{} \  J\ll\log\log N,$ 
\begin{gather}
N_j=
\left\{
              \begin{array}{ll}
                N^{\frac{1}{2-\epsilon_0}(1-\epsilon_0)^{1-j}}, & \hbox{если $-1-J\le j\le1$;} \\
                N^{1-\frac{1}{2-\epsilon_0}(1-\epsilon_0)^{j}}, & \hbox{если $0\le j\le J$,}
              \end{array}
\right.
\end{gather}
имеющая такое свойство    \cite[лемма 9.1]{FK3}: при $-J-1\le m\le J-1$ выполнены неравенства
\begin{equation}\label{ma}
N_{j-J-1}\ge N_{j-J}^{1-\epsilon_0},{}\ {}\  \frac{N}{N_{j-J}}\ge\left( \frac{N}{N_{j-J-1}}\right) ^{1-\epsilon_0}.
\end{equation}
Также в \cite{FK3} по последовательности $\{N_j\}$  было построено специальное множество матриц --- ансамбль  $\Omega_{N}\subseteq \left\{g\in\Gamma_{\A}\,\Bigl|\,\|g\|\le1,02~ N\right\}$ (терминология и существо дела --- из~\cite{BK}),
 \begin{equation}
\label{normi}
\Omega_N=\Xi_1\Xi_2\ldots\Xi_{2J}\Xi_{2J+1},
\end{equation}
где $J$  ---  из (\ref{hgfj}), а все множества $\Xi_j$  ---  предансамбли. 

Отметим, что для доказательства неравенств (\ref{lld}) и (\ref{cdd}) Хуанг в работе  \cite{Huang}  использовал несколько иной вариант ансамбля $\Omega_N$, который отличается, в основном, тем, что среди множеств-множителей в правой части равенства (\ref{normi}) присутствует некоторое специальное множество матриц $\chi$ (впервые построенное в \cite{BK} и перенесенное Хуангом из начала последовательности множеств, сходной с (\ref{normi}), в ``середину'' той же последовательности). Строго говоря, аналогичное видоизменение ансамбля требуется и в настоящей работе --- для тех же целей. Но поскольку такой подход приводит к существенному усложнению выкладок, будем для краткости всегда полагать, что ансамбль задан равенством  (\ref{normi}).

Напомним основные свойства ансамбля.

\begin{Le}\label{6.1}
 \cite[лемма 11.2]{FK3}.    Для любого набора матриц
$$\xi_1\in\Xi_1,\,\xi_2\in\Xi_2,\ldots,\xi_{2J+1}\in\Xi_{2J+1},$$
для любого $j$ из интервала $1\le j\le 2J+1$ выполнено неравенство
\begin{equation}
\label{normp}
\frac{1}{70A^2} N_{j-J}\le \left\|\xi_1\xi_2\ldots\xi_j\right\|\le 1,01 N_{j-J},
\end{equation}
а для каждого  $j$ из интервала $0\le j\le 2J$ имеет место аналогичное неравенство
\begin{equation}
\label{normee}
\frac{1}{150A^2} \frac{N}{N_{j-J}}\le\|\xi_{j+1}\xi_{j+2}\ldots\xi_{2J+1}\|\le  73A^2\frac{N}{N_{j-J}}.
 \end{equation}
\end{Le}

\begin{Sl}\label{s6.1} Если числа $j$ и $U$ связаны соотношением 
\begin{equation}
\label{norhhm}
N_{j-J-1}\le U <N_{j-J},
 \end{equation}
то в обозначениях предыдущей леммы выполнены неравенства
\begin{equation}
\label{normii}
 \frac{1}{70A^2} U\le\|\xi_1\xi_2\ldots\xi_j\|\le 1,01 U^{1+2\epsilon_0},
\end{equation}
\begin{equation}
\frac{1}{150A^2}\max\left \{\frac{N}{U^{1+2\epsilon_0}}, \left( \frac{N}{U}\right)^{1-\epsilon_0} \right\} \le\|\xi_{j+1}\xi_{j+2}\ldots\xi_{2J+1}\|\le 73A^2\frac{N}{U}.
\label{normio}
 \end{equation}
 \end{Sl}
 
 Доказательство. Из неравенств (\ref{ma}) и (\ref{norhhm}) следует двусторонняя оценка числа $N_{j-J}$:
 $$U\le N_{j-J}\le\min\left \{U^{1+2\epsilon_0},{}\  U^{1-\epsilon_0}N^{\epsilon_0}  \right\}.$$
Подставляя ее в (\ref{normp}) и (\ref{normee}), получаем (\ref{normii}) и (\ref{normio}). Следствие доказано.

Пусть $M^{(1)}>1$ --- некоторое действительное число. Определим числа
\begin{equation}
\label{norkkm}
 H=H_1 \left(M^{(1)}\right)=1,01  \left  (M^{(1)}\right)^{1+2\epsilon_0},
\end{equation}
\begin{equation}
\label{lorm}
Q_1=\exp{ \left  (A^4\epsilon_0^{-5}\right)} 
 \end{equation}
 и интервал
$$I_N= \left  [Q_1,N\right].$$
    
  \begin{Th}\label{t6.1}
 (см.  \cite[теорема 4.2]{FK4}). Пусть числа $$M^{(1)}\in I_N, {}\  \ M^{(2)},M^{(4)}\in I_N\bigcup\{1\}$$ удовлетворяют неравенству
\begin{equation}
\label{yorm}
M^{(1)}M^{(2)}M^{(4)}\le N.
  \end{equation}
  Тогда найдутся целые числа --- индексы $j_1,j_2,j_3,$ такие что 
\begin{equation}
\label{nort}
1=j_0\le j_1\le j_2\le j_3\le j_4=2J+2
  \end{equation}
 и  для $k=1,2,3,4$ в обозначениях
 \begin{equation}
\label{nortm}
 \Omega^{(k)}=\Xi_{j_{k-1}}\Xi_{j_{k-1}+1}\Xi_{j_{k-1}+2}\ldots\Xi_{j_{k}-1}
 \end{equation}
  имеют место соотношения (\ref{1norm}) --- (\ref{2norm}) (ниже), а для любых четырех матриц
\begin{equation}
\label{normyy}
 g_1\in\Omega^{(1)},{} \ g_2\in\Omega^{(2)},{} \ g_3\in\Omega^{(3)}, {} \ g_4\in\Omega^{(4)}
 \end{equation}
 выполняются неравенства (\ref{4norm}) --- (\ref{8norm}) (ниже); в частности, если $M^{(2)}=1$ или $M^{(4)}~=~1,$ то $\Omega^{(2)}=\{E\}$ или $\Omega^{(4)}=\{E\}$, соответственно.
 \end{Th}

  Где:  
 
 --- произведение в (\ref{nortm}) по пустому множеству индексов (при $j_{k-1}=j_{k}$) считается равным $\{E\},{}\  E$ --- единичная матрица размера $2\times2$;
\begin{equation}
\label{1norm}
 \Omega_N=\Omega^{(1)}\Omega^{(2)}\Omega^{(3)}\Omega^{(4)},
 \end{equation}

\begin{equation}
\label{2norm}
 \left  |\Omega^{(1)}\Omega^{(2)}\right |>>{ \left  (M^{(1)}M^{(2)}\right)}^{2\delta_{\A}-\epsilon_0},{}\  
 \left   |\Omega^{(4)}\right|>>{ \left  (M^{(4)}\right)}^{2\delta_{\A}-2\epsilon_0},
 \end{equation}
 
\begin{equation}
\label{4norm}
\frac{M^{(1)}}{70A^2} \le\|g_1\|\le H,
 \end{equation}
 
\begin{equation}
\label{5norm}
\frac{N}{150A^2H} \le\|g_2g_3g_4\|\le 73A^2\frac{N}{M^{(1)}},
 \end{equation}

\begin{equation}
\label{6norm}
 \|g_2\|\le 73A^2   M^{(2)}{ \left  (M^{(1)}M^{(2)}\right)}^{2\epsilon_0},
 \end{equation}

\begin{equation}
\label{7norm}
\frac{{ \left  (M^{(4)}\right)}^{1-\epsilon_0} }{150A^2} \le\|g_4\|\le 73A^2M^{(4)},
 \end{equation}

\begin{equation}
\label{8norm}
 \frac{M^{(1)}M^{(2)}}{70A^2} \le\|g_1g_2\|\le 1,01{ \left  (M^{(1)}M^{(2)}\right)}^{1+2\epsilon_0}.
\end{equation}

Доказательство. Положим 
$$U^{(1)}=M^{(1)},{}\  U^{(2)}=M^{(1)}M^{(2)},{}\   U^{(3)}=\frac{N}{M^{(4)}},{}\   U^{(4)}=N. $$
Тогда, ввиду (\ref{yorm}), имеет место цепочка неравенств:
\begin{equation}
\label{9norm}
U^{(1)}\le U^{(2)}\le U^{(3)}\le U^{(4)}.
  \end{equation}
  Поэтому, если для $k=1,2,3$ числа $j_k$ определить неравенствами 
  $$N_{j_k-J-1}\le U^{(k)} <N_{j_k-J},$$
  то, ввиду (\ref{9norm}) и монотонности последовательности $\{N_j\}$, имеют место неравенства  (\ref{nort}). Следовательно, из (\ref{normi}) и (\ref{nortm}) можно вывести равенство (\ref{1norm}). Далее, условия следствия \ref{s6.1} выполнены  при $k=1,2,3$ для $U=U^{(k)}$ и $j=j_k$, для которых, поэтому, имеют место оценки (\ref{normii}) и (\ref{normio}). Первая из них при подстановках $U=U^{(1)}=M^{(1)}$ и $U~=~U^{(2)}=M^{(1)}M^{(2)}$ дает неравенства (\ref{4norm}) и (\ref{8norm}),  соответственно. Аналогично, вторая из оценок следствия \ref{s6.1} дает неравенства (\ref{5norm}) и (\ref{7norm}) при подстановках $U=U^{(1)}=M^{(1)}$ и $U=U^{(3)}=\frac{N}{M^{(4)}}$, соответственно. Комбинация доказанных оценок (\ref{4norm}) и (\ref{8norm}) дает неравенство (\ref{6norm}), так как, в виду нижней оценки в (\ref{continuant inequality}) и равенства (\ref{matrix-normm}), выполняется оценка 
  $$\|g_2\|\le\frac{\|g_1g_2\|}{\|g_1\|}.$$
  Наконец, неравенства (\ref{2norm})  выводятся с помощью теоремы \cite[теорема 12.4]{FK3}, исходя из нижних оценок в доказанных неравенствах (\ref{7norm}) и (\ref{8norm}). Теорема доказана.
  
  В условиях последней теоремы положим: 
\begin{equation}
\label{10norm}
\Omega=\Omega^{(2)}\Omega^{(3)}\Omega^{(4)}=\Xi_{j_{1}}\Xi_{j_{1}+1}\ldots\Xi_{2J+1}.
 \end{equation}
Договоримся использовать следующие обозначения: если $g$ --- матрица размера $2\times2$, то $\widetilde{g}$ --- вектор-столбец, получающийся из матрицы $g$ домножением на $\left (0\atop 1\right)$:

$$\widetilde{g}=g\begin{pmatrix}
0   \\
1
\end{pmatrix}.$$
     Кроме того, если $X$ --- некоторое множество матриц $g$   размера $2\times2$, то $\widetilde{X}$ --- множество соответствующих им вектор-столбцов $\widetilde{g}$.
  Далее, выберем произвольный элемент $g_3\in\Omega^{(3)}$ и положим 
\begin{gather}
\label{ghjk}
\Omega(g_3)=\Omega^{(2)}\{ g_3\}\Omega^{(4)}
=\left\{{g}=g_2g_3g_4 \Bigl | \ g_2\in\Omega^{(2)},{} \ g_4\in\Omega^{(4)} \right\}.
\end{gather}
Тогда, если $\widetilde{g}\in \widetilde{\Omega}(g_3)$, то $\widetilde{g}$ --- вектор-столбец. Для двух произвольных векторов из $\widetilde{\Omega}(g_3)$ введем обозначения:
\begin{equation}
\label{11norm}
 \widetilde{g}^{(1)}={g_2}^{(1)}g_3{g_4}^{(1)}
\begin{pmatrix}
0   \\
1
\end{pmatrix}=\begin{pmatrix}
 x_1  \\
x_2
\end{pmatrix},{}\ 
\widetilde{g}^{(2)}={g_2}^{(2)}g_3{g_4}^{(2)}
\begin{pmatrix}
0   \\
1
\end{pmatrix}=\begin{pmatrix}
y_1 \\
y_2
\end{pmatrix}.
\end{equation}

\begin{Sl}\label{s6.2} Пусть выполнены условия последней теоремы. Тогда, во-первых, имеет  место неравенство
\begin{equation}
\label{12norm}
 \max\limits_{g\in\Omega}\|g\|\le11000A^4\min\limits_{g\in\Omega}\|g\|;
 \end{equation}
 во-вторых, для векторов (\ref{11norm}) для каждого $i\in\{1,2\}$ выполняются неравенства:
\begin{equation}
\label{13norm}
 \max\left \{x_i,y_i \right\}\le22000A^5x_i,
 \end{equation}
\begin{equation}
\label{14norm}
 \frac{N}{2H}\le150A^3x_i.
 \end{equation} \end{Sl}
 
 Доказательство. Согласно лемме \ref{6.1}, в обозначениях (\ref{10norm}) при $j=j_1$ имеет место неравенство (\ref{normee}), в котором верхняя оценка отличается от нижней не более, чем в $11000A^4$ раз. Этим обстоятельством доказаны как неравенство (\ref{12norm}), так и оценка  (\ref{13norm}) при $i=2$ ввиду (\ref{matrix-normm}). Далее, из (\ref{continuant inequality})  следует, что $x_1\le x_2\le 2Ax_1$. Поэтому, ввиду доказанного неравенства (\ref{12norm}), имеет место цепочка оценок:
 $$\max\left   \{x_1,y_1\right\}\le\max \left  \{x_2,y_2\right\}\le11000A^4\min \left  \{x_2,y_2\right\}\le11000A^4x_2\le22000A^5x_1, $$
 так что неравенство (\ref{13norm}) доказано и при $i=1$. Наконец, если в нижней из оценок (\ref{5norm}) положить $g_2=g^{(1)}_2,{}\ g_4=g^{(1)}_4$ и учесть обозначения (\ref{11norm}), то получим: 
 $$\frac{N}{2H}\le75A^2x_2\le150A^3x_1\le150A^3x_i$$
 для каждого $i\in\{1,2\}$, и неравенство (\ref{14norm}) также доказано. Следствие доказано.
 
 \section{Достаточные условия                выполнения оценок (\ref{lc})~---~(\ref{cdd}).}\label{7}
 
 Согласно методу Бургейна --- Конторовича  \cite{BK}, оценка   (\ref{lc}) следует из неравенства 
\begin{equation}
\label{15norm}
 \int^{1}_{0}|S_N( \Theta)|^2d\Theta\ll\frac{\left   |\Omega_N\right|^2}{N},
 \end{equation}
 где 
 \begin{equation}
 S_N( \Theta)=\sum_{g\in\Omega_N}e(\Theta\|g\|)=\sum_{g\in\Omega_N}e\left((0,1) g\begin{pmatrix}
0   \\
1\end{pmatrix}\Theta\right)=\sum_{g\in\Omega_N}e((0,1)\widetilde{g}\Theta)
\end{equation} 
--- тригонометрическая сумма. 

 Применяя теорему Дирихле \cite[лемма 2.1, стр.17]{Von}, для каждого $\Theta\in[0,1)$ найдем целые числа $a$ и $q$  и действительное число $K$, такие что 
\begin{equation}
\label{17norm}
  \Theta=\left\{\frac{a}{q}+\frac{K}{N}\right\},{}\ (a,q)=1, {}\ 0\le a<q\le\frac{N^{\frac{1}{2}}}{Q_1},{}\ {}\ |K|\le\frac{Q_1N^{\frac{1}{2}}}{q}
 \end{equation}
(где фигурные скобки обозначают дробную долю, а величина $Q_1$ была определена в  (\ref{lorm})), при чем равенство $a=0$ возможно только при $q=1$. Далее, представим число $K$ в виде 

\begin{equation} 
\label{19norm}
K=\frac{1}{2}l+\lambda,{}\  \lambda\in \left (-\frac{1}{4},\frac{1}{4}\right ],
  \end{equation}
где $l$ --- какое-либо целое число, тогда $|l|\le2|K|+1$, поэтому, ввиду  \ref{17norm}, имеет место оценка
\begin{equation} 
\label{20norm}
 |l|\le\frac{3}{q}Q_1N^{\frac{1}{2}}.
 \end{equation}
Для произвольного неотрицательного параметра $Q=Q(N)$ определим величину $\Sigma_N(Q)$:
\begin{equation} 
\label{21norm}
 \Sigma_N(Q)=\int^{\frac{1}{4}}_{-\frac{1}{4}}{}\ 
 \mathop{{\sum}^*}_{0\le a<q\le\frac{1}{Q_1}N^{\frac{1}{2}}}{}\ 
 \sum_{|l|\le\frac{3}{q}Q_1N^{\frac{1}{2}}    }{}\ 
 \1_{\{\max\{q,|l|\}\ge Q\}}\left   |S_N\left   ( \frac{a}{q}+\frac{l}{2N}+\frac{\lambda}{N}\right)\right|^2d\lambda,
  \end{equation}
где сумма $\mathop{{\sum}^*}$ берется по парам взаимно простых чисел $a$ и $q$. 

\begin{Le}\label{7.1}  Неравенство (\ref{lc}) имеет место, если найдется хотя бы одно действительное число $\epsilon_0~\in~\left  (0,\frac{1}{2500}\right)$, такое что выполнено неравенство 
\begin{equation} \vec{}
\label{22norm}
 \Sigma_N(0)\ll\left   |\Omega_N\right|^2.
 \end{equation}
 \end{Le}
 
 Доказательство. Положим $Q=0$ и оценим сверху интеграл из (\ref{15norm})  интегралом из (\ref{21norm}), деленным на $N$. Для этого достаточно представить $\Theta$ в виде (\ref{17norm}), где число $K$ взято из (\ref{19norm}), разбить отрезок интегрирования в  (\ref{15norm}) на ряд отрезков, соответствующих значениям параметров $a,q,l,$ и в получившихся интегралах сделать линейную замену $\lambda=N \left( \Theta -\frac{a}{q}-\frac{l}{2N}\right) $. Лемма доказана.
 
 \begin{Zam}\label{z7.1}  Согласно результатам работы  \cite{Huang}, неравенства (\ref{lld}) и (\ref{cdd}) имеют место, если найдется хотя бы одно действительное число $\epsilon_0~\in~\left  (0,\frac{1}{2500}\right)$, такое что для некоторой константы $C=C(\A)>0$ для любой сколь угодно медленно стремящейся к бесконечности величины $Q=Q(N)$ выполнено неравенство
\begin{equation} 
\label{23norm}
\Sigma_N(Q)\ll\left   |\Omega_N\right|^2Q^{-C}.
  \end{equation}\end{Zam}
  
  Напомним, что $Q_1$ было определено в (\ref{lorm}), положим $Q_0=0$    и определим последовательность $\{Q_j\}$ для $j$ от нуля до бесконечности:
 \begin{equation} 
\label{24norm}
  \left  \{Q_j\right\}^\infty_{j=0}= \left  \{0,Q_1,Q^2_1,Q^3_1,\ldots,Q^j_1,\ldots\right\}.
 \end{equation}
 
 При фиксированном $\lambda\in\left (-\frac{1}{4},\frac{1}{4}\right ]$ для целых $\alpha,\beta\ge0$ положим:
 \begin{equation} 
\label{25norm}
P^{(\lambda)}_{\alpha,\beta}=\left \{\Theta=\frac{a}{q}+\frac{l}{2N}+\frac{\lambda}{N}\Bigl|{} \ (a,q)=1, {} \ 0\le a< q,{}\ \mbox{выполнены}{}\  (\ref{26norm}){}\  \mbox{и}{}\  (\ref{27norm}) \right\},
\end{equation}
где
\begin{equation} 
\label{26norm}
Q_ \alpha\le q\le \min\left\{  \frac{1}{Q_1}N^{\frac{1}{2}},Q_{\alpha+1} \right\}, 
\end{equation}

\begin{equation} 
\label{27norm}
Q_ \beta\le         \left|l \right|         \le \min\left\{  \frac{3}{q}Q_1N^{\frac{1}{2}},Q_{\beta+1} \right\}.
  \end{equation}
  
  \begin{Zam}\label{z7.2}  Заметим, что для чисел $Q_{\alpha+1}$,  $Q_{\beta+1}$, при $\alpha,\beta\ge0$ для всех непустых множеств $P^{(\lambda)}_{\alpha,\beta}$ выполняется неравенство
\begin{equation} 
\label{28norm}
 Q_{\alpha+1}Q_{\beta+1}\le 3Q_3N^{\frac{1}{2}}.
 \end{equation}
 Действительно, пусть  какое-либо $\Theta$ из $P^{(\lambda)}_{\alpha,\beta}$ имеет параметры $a$, $q$ и $l$. Тогда из (\ref{26norm}) и (\ref{27norm}) следует, что  
 $$Q_{\alpha+1}\le Q_1q,{} \   Q_{\beta+1}\le Q_1\left|l \right|  $$
по определению последовательности (\ref{24norm}). Отсюда и из неравенства (\ref{20norm}) (имеющего место ввиду (\ref{27norm})) следует утверждение замечания.\end{Zam}

Определим действительное число $\gamma$ равенством
 \begin{equation} 
\label{29norm}
\gamma=1-\delta_{\A}
  \end{equation}
  и для всякого действительного $x$ положим
\begin{equation} 
\label{30norm}
\overline{x}=\max \{1, |x|\}.
 \end{equation}
 
 \begin{Le}\label{7.2} Пусть найдутся константа $c_1=c_1(\A,\epsilon_0)>0$ и абсолютная постоянная величина $c_2>0$, такие что для любого $ \lambda\in \left (-\frac{1}{4},\frac{1}{4}\right ]$, для любых целых  $\alpha,\beta\ge0$, при любом $\epsilon_0~\in~\left  (0,\frac{1}{2500}\right)$ выполняется неравенство 
\begin{equation} 
\label{31norm}
 {\sum_{\Theta\in P^{(\lambda)}_{\alpha,\beta}} |S_N( \Theta)|^2}\le
  c_1 \frac{\left   |\Omega_N\right|^2Q^2_{\alpha+1}Q_{\beta+1}}
  {\overline{Q}^3_{\alpha}\overline{Q}^2_{\beta}}
    {\left(Q_{\alpha+1}Q_{\beta+1}\right)}^{ c_2\gamma+100\epsilon_0}.
\end{equation}
  Тогда из оценки
  \begin{equation} 
\label{32norm}
\delta_{\A}>1-\frac{1}{c_2}
  \end{equation}
следуют неравенства (\ref{lc})~---~(\ref{cdd}).\end{Le}

Доказательство. Предположим, что неравенства (\ref{31norm}) и (\ref{32norm}) выполнены, и докажем оценки  (\ref{22norm}) и  (\ref{23norm}) для некоторых величин $\epsilon_0~\in~\left  (0,\frac{1}{2500}\right)$ и $C=C(\A)>0$. Оценим сумму $\Sigma_N(Q)$, определенную в (\ref{21norm}), с помошью введенного в (\ref{25norm}) множества $P^{(\lambda)}_{\alpha,\beta}$:
\begin{equation} 
\label{34norm}
\Sigma_N(Q)\le \int^{\frac{1}{4}}_{-\frac{1}{4}} \sum_{\alpha,\beta\ge0}\1_{\left\{Q_{\alpha+1}\ge Q\mbox{или } Q_{\beta+1} \ge Q\right\}}  
 \sum_{\Theta \in P^{(\lambda)}_{\alpha,\beta}}
 \left   |S_N(\Theta)\right|^2d\lambda, 
 \end{equation}
 где сумма по $\alpha,\beta$ распространена на целые неотрицательные числа. Подставляя в (\ref{34norm}) равномерную по $\lambda$ оценку  (\ref{31norm}), получаем:
\begin{equation} 
\label{35norm}
 \Sigma_N(Q)\le c_1Q_3\left   |\Omega_N\right|^2\sum_{\alpha,\beta\ge0}\1_{\left\{Q_{\alpha+1}\ge Q\mbox{ или }Q_{\beta+1}\ge Q\right\}}\left( Q_{\alpha+1}Q_{\beta+1}\right)^{-1+c_2\gamma+100\epsilon_0}.
 \end{equation}
Определим теперь числа $\epsilon_0$ и $C=C(\A)$ равенствами
$$\epsilon_0=\min\left\{\frac{1}{200}(1-c_2\gamma),{}\ \frac{1}{2600}\right\},$$
\begin{equation} 
\label{37norm}
 C=C(\A)=1-c_2\gamma-100\epsilon_0.
 \end{equation}
Тогда, ввиду соотношения (\ref{32norm}),  выполняются неравенства $\epsilon_0>0$ и $C>0$. Определим еще две величины: целое число $\xi=\xi(Q)\ge0 $ и
 для произвольных $\alpha_0,\beta_0\ge0$ --- сумму $\Sigma_N(\alpha_0,\beta_0)$, соответственно, из условий:
\begin{equation} 
\label{40norm}
Q_\xi\le Q<Q_{\xi+1},
  \end{equation}
 $$ \Sigma_N(\alpha_0,\beta_0)=\sum_{\alpha\ge\alpha_0}\sum_{\beta\ge\beta_0}\left( Q_{\alpha+1}Q_{\beta+1}\right)^{-C}.$$
Тогда, ввиду соотношений (\ref{35norm}) и (\ref{37norm}), выполняется оценка
\begin{equation} 
\label{39norm}
\Sigma_N(Q)\ll \left   |\Omega_N\right|^2\left( \Sigma_N(0,\xi)+ \Sigma_N(\xi,0) \right).
  \end{equation}
Заметим, что из неравенства (\ref{40norm}) следует оценка $Q_{\xi+1}\ge\overline{Q}$. В частности, если $\Sigma_N(\alpha_0,\beta_0)$ --- любое из двух слагаемых в скобках правой части неравенства (\ref{39norm}), то $Q_{\alpha_0+1}\ge\overline{Q}$ или $ Q_{\beta_0+1}\ge\overline{Q}$. Поэтому из оценки произведения сумм сходящихся геометрических прогрессий
$$ \Sigma_N(\alpha_0,\beta_0)= \left( \sum_{\alpha\ge\alpha_0}\left( Q_{\alpha+1}\right)^{-C}\right) \left(\sum_{\beta\ge\beta_0}\left(Q_{\beta+1}\right)^{-C}  \right)\ll\left( Q_{\alpha_0+1}Q_{\beta_0+1}\right)^{-C} $$
следуют неравенства
$$\Sigma_N(0,\xi)\ll\left(  \overline{Q} \right)^{-C},{}\ \Sigma_N(\xi,0)\ll\left(  \overline{Q} \right)^{-C} .$$
Подставляя эти оценки в неравенство (\ref{39norm}), получаем соотношения (\ref{22norm}) и (\ref{23norm}). Теперь утверждение настоящей леммы леммы следует из леммы \ref{7.1} и замечания \ref{z7.1}. Лемма доказана.

Фиксируем константу
 \begin{equation}
\label{43norpm}
T_1=7Q_7
\end{equation}
(где $Q_7$ --- элемент последовательности (\ref{24norm})) и для любого целого $\kappa\in\left[0,{}\ T_1-1\right]$ положим:
\begin{equation}
\label{43norm}
P^{(\lambda)}_{\alpha,\beta}(\kappa)=\left \{\Theta\in P^{(\lambda)}_{\alpha,\beta}\Bigl|{} \ l\equiv\kappa\pmod {T_1}\right\}.
  \end{equation}

Всюду далее $\alpha,\beta\ge0, {}\ 0\le\kappa<T_1$ --- произвольные целые числа, {}\ $\lambda\in \left (-\frac{1}{4},\frac{1}{4}\right ]$ --- произвольное действительное. Через $Z$ всюду далее обозначается произвольное непустое подмножество конечного множества $P^{(\lambda)}_{\alpha,\beta}(\kappa)$. Положим
\begin{equation} 
\label{44norm}
\sigma_{N,Z}=\sum_{\Theta\in Z}\left|S_N(\Theta)\right|.
\end{equation}
\begin{Th}\label{t7.1} Если найдутся константа $c_1=c_1(\A, \epsilon_0)>0$ и абсолютная постоянная величина $c_2~>~0$, такие что при любом $\epsilon_0~\in~\left  (0,\frac{1}{2500}\right)$ выполняется неравенство  
\begin{equation} 
\label{45norm}
 \sigma_{N,Z}\le
  (c_1)^{\frac{1}{2}} \frac{\left   |\Omega_N\right||Z|}
  {(\overline{Q}_{\alpha}\overline{Q}_{\beta}|Z|)^{\frac{1}{2}}}
    {\left(Q_{\alpha+1}Q_{\beta+1}\right)}^{ \frac{1}{2}\left (c_2\gamma+{90}\epsilon_0\right)},
\end{equation}
то из оценки (\ref{32norm}) следуют неравенства (\ref{lc})~---~(\ref{cdd}).\end{Th} 

Доказательство. Достаточно из неравенства (\ref{45norm}) вывести оценку вида (\ref{31norm}) с тем же самым значением $c_2$, тогда утверждение теоремы  будет следовать из леммы 7.2. Вывод неравенства (\ref{31norm}) легко получается из обобщенной леммы Конягина \cite[лемма 13.5]{FK3}: для всякой функции $f:W\rightarrow \rr_{+}$, где $W$ --- конечное множество, имеет место неравенство
\begin{equation} 
\label{46norm}
\sum_{\Theta\in W}f^2(\Theta)\ll_{\epsilon} \left|W \right|^{\epsilon}\max_{Z\subseteq W}\left(\frac{1}{|Z|}\left(\sum_{\Theta\in Z}f(\Theta)  \right)^2 \right),
  \end{equation}
  где максимум берется по всем непустым подмножествам. Для завершения доказательства остается в (\ref{46norm}) положить
$$f(\Theta)=\left|S_N(\Theta)\right|,{}\ W=P^{(\lambda)}_{\alpha,\beta}(\kappa)$$ 
и воспользоваться оценкой  (\ref{45norm}) и неравенством
$${\sum_{\Theta\in P^{(\lambda)}_{\alpha,\beta}} |S_N( \Theta)|^2}\le T_1\max_{0\le\kappa\le T_1-1}\sum_{\Theta\in P^{(\lambda)}_{\alpha,\beta}(\kappa)} |S_N( \Theta)|^2.$$
Теорема доказана.

Цель всех последующих рассуждений будет состоять в доказательстве неравенства(\ref{45norm}) при $c_2=5$.

\section{Общие оценки тригонометрической суммы по ансамблю.}

Всюду в данном параграфе считаем, что числа $M^{(1)},M^{(2)},M^{(4)}$, удовлетворяющие условиям теоремы 6.1, уже каким-либо образом выбраны, так что имеет место разложение ансамбля (\ref{1norm}) (со свойствами (\ref{2norm}) --- (\ref{8norm})). Именно вид такого разложения отличает данный параграф от аналогичного ему \cite[параграф 13]{FK3} с тем же названием. Ввиду почти дословного соответствия формулировок, доказательство первых двух лемм настоящего параграфа не приводится.

Далее используются обозначения: (\ref{norkkm}), (\ref{ghjk}) и (\ref{11norm}) из параграфа \ref{6} , (\ref{25norm}) и (\ref{44norm}) из парагрфа \ref{7}, а также обозначения из \cite{FK3}:
$$S(x)=3\left(   \frac{2\sin\left(\frac{1}{2}\varPi x \right) }{\varPi x}\right)^2, {}\ \$(z)=S(z_1)S(z_2)$$
для  $z=(z_1,z_2)\in \rr^2$ при $z\neq (0,0)$, и $\$(0,0)=3$ --- по непрерывности.
 
 \begin{Le}\label{8.1}  (см. \cite[лемма 13.1]{FK3}). При некоторых комплексных коэффициентах $\xi (\Theta)$ (где $\Theta\in Z$), по модулю равных единице, имеет место оценка
$$\sigma_{N,Z}\ll\left   |\Omega^{(1)}\right|^{\frac{1}{2}}
\sum_{g_3\in \Omega^{(3)}}\left(\sum_{g_1\in \mathbb{Z}^2}\$\left(\frac{g_1}{H} \right)\left|\sum_{\Theta\in Z}\xi (\Theta)\sum_{\widetilde{g}\in \widetilde{\Omega}(g_3)}
e(g_1\widetilde{g}\Theta)   \right|^2\right)^{\frac{1}{2}}.$$\end{Le}

Отсюда легко получается следующая лемма.

\begin{Le}\label{8.2}  (см. \cite[лемма 13.2]{FK3}). Имеет место оценка
\begin{equation} 
\label{49norm}
 \sigma_{N,Z}\ll H\left   |\Omega^{(1)}\right|^{\frac{1}{2}}
\sum_{g_3\in \Omega^{(3)}}
\left(\sum_{{g^{(1)},g^{(2)}\in \widetilde{\Omega}(g_3)}\atop{\Theta^{(1)},\Theta^{(2)}\in Z}}
\1_
{\left \{
 \left\|   
 \widetilde{g}^{(1)}\Theta^{(1)}-
\widetilde{g}^{(2)}\Theta^{(2)}  \right\|    
_{1,2}\le\frac{1}{2H}\right\}}\right)
^{\frac{1}{2}},
\end{equation}
где $\left\|z\right\|_{1,2}=\max \left\{   \left\|z_1\right\|,\left\|z_2\right\|\right\}$ для $z=(z_1,z_2)^t\in \rr^2$.\end{Le}

Напомним, что всегда $Z\subseteq P^{(\lambda)}_{\alpha,\beta}(\kappa)$ для соответствующих значений параметров. Числа $\Theta_1,\Theta_2\in P^{(\lambda)}_{\alpha,\beta}(\kappa)$ будем далее записывать следующим образом:
\begin{equation} 
\label{50norm}
\Theta^{(1)}=\frac{a^{(1)}}{q^{(1)}}+\frac{l^{(1)}}{2N}+\frac{\lambda}{N}, {}\ {}\ \Theta^{(2)}=\frac{a^{(2)}}{q^{(2)}}+\frac{l^{(2)}}{2N}+\frac{\lambda}{N}.
\end{equation}
Положим:
$$\NN(g_3)=\left\{    \left(g^{(1)},g^{(2)},\Theta^{(1)},\Theta^{(2)}\right) \in   \left(\widetilde{\Omega}(g_3)\right)^2\times Z^2 \Bigl|{} \  \mbox{(\ref{52norm})}{}\  \mbox{и}{}\ \mbox{(\ref{53norm})}{}\ \mbox{выполнены для}{}\ i=1,2 \right\},$$
где (в обозначениях из (\ref{11norm})):
\begin{equation} 
\label{52norm}
\left\|  \frac{x_ia^{(1)}}{q^{(1)}}-\frac{y_ia^{(2)}}{q^{(2)}}\right\|\le\min
\left\{
              \begin{array}{ll}
                \frac{1}{N}150A^3x_i+\left\|  \frac{1}{2N}\left(x_il^{(1)}-y_il^{(2)}\right)\right\|+
                \left\| \frac{\lambda}{N}\left(  x_i-y_i\right)\right\|,                \\
              \frac{1}{N} \left(9A \right) ^5x_iQ_{\beta+1},\\
                \frac{1}{M^{(1)}}74A^2Q_{\beta+1},
              \end{array}
\right.
\end{equation}

\begin{equation} 
\label{53norm}
\left| x_il^{(1)}-y_il^{(2)} \right|\le\left(9A \right) ^5x_i+2N\left\|  x_i\frac{a^{(1)}}{q^{(1)}}-y_i\frac{a^{(2)}}{q^{(2)}}\right\|.
\end{equation}

\begin{Le}\label{8.3}  (см. \cite[лемма 13.3]{FK3}). Если выполнено неравенство
\begin{equation} 
\label{54norm}
M^{(1)}\ge150A^2Q_{\beta+1},
\end{equation}
то имеет место оценка
\begin{equation} 
\label{55norm}
\sigma_{N,Z}\ll H\left   |\Omega^{(1)}\right|^{\frac{1}{2}}
\sum_{g_3\in \Omega^{(3)}}\left| \NN(g_3)   \right|^{\frac{1}{2}}.
\end{equation}\end{Le}

Доказательство. Будем выводить оценки (\ref{52norm}) и (\ref{53norm}) из неравенства в фигурных скобках в (\ref{49norm}): тогда утверждение настоящей леммы будет следовать из предыдущей. 

Применение верхней оценки в (\ref{27norm}) или неравенства $\left| \lambda \right|\le\frac{1}{4}$ дает, соответственно, неравенства:
$$\left|\frac{ x_il^{(1)}}{2N}-\frac{ y_il^{(2)}}{2N} \right|\le\frac{Q_{\beta+1}}{N}\max \left\{  x_i,y_i \right\},$$

$$\left|\frac{ x_i\lambda}{N}-\frac{ y_i\lambda}{N} \right|\le\frac{2\left| \lambda \right|}{N}\max \left\{  x_i,y_i \right\}
\le\frac{1}{2N}\max \left\{  x_i,y_i \right\}.$$
Ввиду обозначений, введенных равенствами (\ref{10norm}) и (\ref{ghjk}), имеет место включение $\Omega(g_3)\subseteq\Omega$. Следовательно, к максимуму из двух последних неравенств можно применить любую из двух верхних оценок: ту, которая в (\ref{5norm}), или (\ref{13norm}). Их применение дает:
\begin{equation} 
\label{58norm}
\left|\frac{ x_il^{(1)}}{2N}-\frac{ y_il^{(2)}}{2N} \right|\le\min \left\{  \frac{73A^2Q_{\beta+1}}{M^{(1)}},{}\ 
\frac{1}{N}22000A^5Q_{\beta+1}x_i \right\},
\end{equation}
\begin{equation} 
\label{59norm}
\left\|\frac{ x_i\lambda}{N}-\frac{ y_i\lambda}{N} \right\|      \le\min \left\{  \frac{73A^2}{2M^{(1)}},{}\ 
\frac{1}{N}11000A^5x_i \right\}.
\end{equation}
Из неравенства (\ref{54norm}) следует, что первый элемент минимума в (\ref{58norm}) меньше, чем $\frac{1}{2}$, поэтому
\begin{equation} 
\label{60norm}
\left\|\frac{ x_il^{(1)}}{2N}-\frac{ y_il^{(2)}}{2N} \right\|=\left|\frac{ x_il^{(1)}}{2N}-\frac{ y_il^{(2)}}{2N} \right|.
\end{equation}
Применяя дважды неравенство треугольника, получаем:
\begin{equation} 
\label{61norm}
\left\|  \frac{x_ia^{(1)}}{q^{(1)}}-\frac{y_ia^{(2)}}{q^{(2)}}\right\|\le
\left\|\frac{ x_il^{(1)}}{2N}-\frac{ y_il^{(2)}}{2N} \right\|+
\left\|\frac{ x_i\lambda}{N}-\frac{ y_i\lambda}{N} \right\| +
\left\|\ x_i\Theta^{(1)}-y_i\Theta^{(2)}\right\|, 
\end{equation}
\begin{equation} 
\label{62norm}
\left\|\frac{ x_il^{(1)}}{2N}-\frac{ y_il^{(2)}}{2N} \right\|\le
\left\|  \frac{x_ia^{(1)}}{q^{(1)}}-\frac{y_ia^{(2)}}{q^{(2)}}\right\|+
\left\|\frac{ x_i\lambda}{N}-\frac{ y_i\lambda}{N} \right\| +
\left\|\ x_i\Theta^{(1)}-y_i\Theta^{(2)}\right\|.
\end{equation}

Заметим, что из (\ref{61norm}) сразу получается верхняя строчка из минимума в (\ref{52norm}), если третье слагаемое в правой части (\ref{61norm}) оценить, как в неравенстве (\ref{49norm}):
\begin{equation} 
\label{63norm}
\left\|\ x_i\Theta^{(1)}-y_i\Theta^{(2)}\right\|\le\frac{1}{2H} \le\frac{1}{N}150A^3x_i
\end{equation}
ввиду  оценки (\ref{14norm}). Если теперь остальные слагаемые из правой части (\ref{61norm}) оценить вторыми элементами минимумов в неравенствах (\ref{58norm}) и (\ref{59norm}), то получим вторую строку минимума в (\ref{52norm}). Далее, третья строчка минимума в (\ref{52norm}) получается применением в неравенстве (\ref{61norm}) первых элементов минимумов из неравенств (\ref{58norm}) и (\ref{59norm}). Теперь неравенство (\ref{52norm}) полностью доказано.

Аналогично, оценка  правой части (\ref{62norm}) производится на основании второго элемента минимума в (\ref{59norm}), а также неравенства (\ref{63norm}). К результату применяется равенство (\ref{60norm}), приводящее к оценке (\ref{53norm}). Лемма доказана.

Обозначим:
$$\M(g_3)=\left\{    \left(\widetilde{g}^{(1)},\widetilde{g}^{(2)},\Theta^{(1)},\Theta^{(2)}\right) \in   \left(\widetilde{\Omega}(g_3)\right)^2\times Z^2 \Bigl|{} \  (\ref{65norm}){}\  \mbox{и}{}\ (\ref{66norm}){}\ \mbox{выполнены для}{}\ i=1,2 \right\},$$
где (в обозначениях (\ref{11norm}) и (\ref{50norm}))
\begin{equation} 
\label{65norm}
\left\|  \frac{x_ia^{(1)}}{q^{(1)}}-\frac{y_ia^{(2)}}{q^{(2)}}\right\|=0,
\end{equation}
\begin{equation} 
\label{66norm}
\left| x_il^{(1)}-y_il^{(2)} \right|\le\left(9A \right) ^5x_i.
\end{equation}

Далее через $\q$ обозначается наименьшее общее кратное чисел $q^{(1)}$ и $q^{(2)}$.

\begin{Zam}\label{z8.1}  В работе \cite[пп. 6.2.1]{BK} (см. также \cite[доказательство леммы 14.1]{FK3}) было показано, что условие (\ref{65norm}) равносильно следующему:
\begin{equation} 
\label{67norm}
q^{(1)}=q^{(2)}=\q ,{}\ x_ia^{(1)}\equiv y_ia^{(2)}\pmod{\q},{}\ i=1,2.
\end{equation}
Отсюда следует, что для элементов множества $\M(g_3)$ из равенства $x_1=y_1$ следуют соотношения: $a^{(1)}=a^{(2)}$ --- в виду сравнения в (\ref{67norm}) при $i=1$, а равенство $l^{(1)}=l^{(2)}$ --- в виду неравенства (\ref{66norm}) при $i=1$ и сравнения в (\ref{43norm}).\end{Zam}

\begin{Le}\label{8.4}(см. \cite[лемма 13.4]{FK3}). Пусть число $M^{(1)}$ определено равенством
\begin{equation} 
\label{68norm}
M^{(1)}=150A^2Q^2_{\alpha+1}Q_{\beta+1}.
\end{equation}
Тогда при любом $g_3\in\Omega^{(3)}$ имеет место включение
\begin{equation} 
\label{69norm}
\NN(g_3)\subseteq\M(g_3).
\end{equation}\end{Le}

Доказательство. Покажем, что соотношения (\ref{65norm}) и (\ref{66norm}) в условиях леммы следуют из неравенств (\ref{52norm}) и (\ref{53norm}). Действительно, в виду (\ref{68norm}), третья строка минимума в правой части неравенства (\ref{52norm}) меньше, чем $\frac{1}{q_1q_2}$ --- следовательно, оцениваемая величина равна нулю. Отсюда следует равенство (\ref{65norm}), подстановка которого в  (\ref{53norm}) дает неравенство(\ref{66norm}). Лемма доказана.

Введем следующие обозначения (для переменных из (\ref{11norm})):
\begin{equation} 
\label{70norm}
\mathcal{Y}=x_1y_2-x_2y_1, 
\end{equation}
\begin{equation}
\label{76norm} 
{}\  {}\ Y=Y(M^{(1)})=\frac{75A^2N}{M^{(1)}}.
\end{equation}

\begin{Le}\label{8.5}(см. \cite[доказательство предложения 6.11]{BK}). Пусть число $M^{(1)}$ определено равенством
\begin{equation} 
\label{71norm}
M^{(1)}=120A^2(NQ_{\alpha+1}Q_{\beta+1})^{\frac{1}{2}}.
\end{equation}
Тогда при любом $g_3\in\Omega^{(3)}$ для элементов множества $\NN(g_3)$ имеют место соотношения

\begin{equation} 
\label{75norm}
0<x_i<Y,{}\ {}\ 0<y_i<Y,{}\ {}\ i=1,2,
\end{equation}
\begin{equation} 
\label{72norm}
\mathcal{Y}<Y^2,
\end{equation}
\begin{equation} 
\label{73norm}
\frac{2Y^2Q_{\beta+1}}{N}<\frac{1}{Q_{\alpha+1}},
\end{equation}
\begin{equation} 
\label{74norm}
\mathcal{Y}\equiv 0 \pmod{\q}.
\end{equation}\end{Le}

Доказательство. Hеравенства (\ref{75norm}) следуют из верхней оценки в (\ref{5norm}), откуда следует неравенство (\ref{72norm}). Далее,  соотношение (\ref{73norm}) проверяется непосредственным применением равенств (\ref{76norm}) и  (\ref{71norm}). Наконец, поскольку
\begin{equation} 
\label{77norm}
 \left\|\mathcal{Y}\frac{a^{(1)}}{q^{(1)}}\right\|\le
y_2 \left\|  \frac{x_1a^{(1)}}{q^{(1)}}-\frac{y_1a^{(2)}}{q^{(2)}}\right\|+
y_1 \left\|\frac{x_2a^{(1)}}{q^{(1)}}-\frac{y_2a^{(2)}}{q^{(2)}}\right\|,
\end{equation}
 то, применяя здесь нижнюю строку минимума в (\ref{52norm})  и оценки (\ref{75norm}), получаем:
\begin{equation} 
\label{78norm}
\left\|\mathcal{Y}\frac{a^{(1)}}{q^{(1)}}\right\|<
2Y\frac{74A^2Q_{\beta+1}}{M^{(1)}}<\frac{2Y^2Q_{\beta+1}}{N}<\frac{1}{Q_{\alpha+1}}<\frac{1}{q^{(1)}}
\end{equation}
ввиду доказанной оценки (\ref{73norm}). Следовательно, $\left\|\mathcal{Y}\frac{a^{(1)}}{q^{(1)}}\right\|=0$, то есть $\mathcal{Y}$ делится на $q^{(1)}$. Аналогично, $\mathcal{Y}$ делится на $q^{(2)}$, откуда следует сравнение (\ref{74norm}). Лемма доказана.

\begin{Le}\label{8.6}  Если число $M^{(1)}$ определено равенством (\ref{71norm}), то при любом $g_3\in\Omega^{(3)}$ имеет место включение (\ref{69norm}).\end{Le}

Доказательство. Как и в доказательстве леммы \ref{8.4}, достаточно лишь вывести равенство (\ref{65norm}) из соотношений (\ref{52norm}), (\ref{53norm}) и условий настоящей леммы.

Пусть для начала $\mathcal{Y}$ не равно нулю (\textbf{случай 1}, см. \cite[пп. 6.2.1]{BK}). Тогда, ввиду соотношений (\ref{72norm}) и (\ref{74norm}), имеет место цепочка неравенств:
\begin{equation} 
\label{79norm}
\q\le\min{\left\{ Q^2_{\alpha+1},|\mathcal{Y}| \right\}}
\le\min{\left\{ Q^2_{\alpha+1},Y^2 \right\}}\le Q_{\alpha+1}Y.
\end{equation}
Отсюда и из (\ref{73norm}) следует, что
\begin{equation} 
\label{80norm}
\frac{YQ_{\beta+1}}{N}\le \frac{Y^2Q_{\alpha+1}Q_{\beta+1}}{\q N}<\frac{1}{\q}.
\end{equation}
Поэтому третья строчка минимума в (\ref{52norm}) меньше, чем $\frac{1}{\q}$. Следовательно, левая часть неравенства (\ref{52norm}) равна нулю, так что имеет место неравенство (\ref{65norm}).

Пусть теперь $\mathcal{Y}$ равно нулю (\textbf{случай 2}, см. \cite[доказательство леммы 14.6]{FK3}). Ввиду взаимной простоты чисел $x_1$ и $x_2$, а также $y_1$ и $y_2$, это означает, что 
\begin{equation} 
\label{81norm}
x_1=y_1,{}\ {}\ x_2=y_2.
\end{equation}
 Поэтому любая из первых двух строчек минимума в (\ref{52norm}) может быть представлена в виде $ F(x_i)$ для $i=1,2$, где $ F$ --- некоторая однородная  линейная функция.  Следовательно, если положить $\xi=\left|\frac{a^{(1)}}{q^{(1)}}-\frac{a^{(2)}}{q^{(2)}}\right|$, то оценка (\ref{52norm}) запишется в виде
$$ \left\|x_i\xi \right\|\le F(x_i).$$
Другими словами, если $n_1$ и $n_2$ --- ближайшие целые числа к величинам $x_1\xi$ и $x_2\xi$, соответственно, то найдутся действительные числа $\Theta_1, \Theta_2$, такие что $|\Theta_1|,|\Theta_2| <\frac{1}{2}$ и
\begin{equation} 
\label{83norm}
n_1=x_1\xi+\Theta_1F(x_1),{}\ {}\ n_2=x_2\xi+\Theta_2F(x_2).
\end{equation}
В частности, если $\ F(x_i)$ --- вторая строка минимума в (\ref{52norm}), то из (\ref{75norm}) и (\ref{83norm}) следует, что
$$|x_1n_2-x_2n_1|=x_1x_2|\Theta_1-\Theta_2|\frac{(9A)^5Q_{\beta+1}}{N}\le
\frac{Y^2(9A)^5Q_{\beta+1}}{N}\le\frac{(9A)^5}{Q_{\alpha+1}}< 1$$
ввиду (\ref{73norm}) (и поскольку $(9A)^5<Q_1\le Q_{\alpha+1}$). Следовательно, $x_1n_2=x_2n_1$. Но числа $x_1$ и $x_2$ взаимно просты, поэтому найдется целое число $k$, такое что $n_2=kx_2$. Отсюда, возвращаясь к случаю произвольного из двух вариантов для $F(x_2)$ во втором из равенств в 
 (\ref{83norm}) и сокращая на $x_2$, получаем:
\begin{equation} 
\label{85norm}
\left\|\xi \right\|\le\left|\xi-k\right|\le\frac{F(x_2)}{x_2}.
\end{equation}

В частности, при $Q_{\alpha+1}<\frac{1}{Q_4}N^{\frac{1}{2}}$ (\textbf{случай 2.1}), если $F(x_i)$ --- вторая строка минимума в (\ref{52norm}), то   из неравенств (\ref{85norm}) и  (\ref{28norm}) получаем:
\begin{equation} 
\label{86norm}
\left\|\xi \right\|\le\frac{(9A)^5Q_{\beta+1}}{N}\le
\frac{3(9A)^5Q_3}{Q_{\alpha+1}N^{\frac{1}{2}}}<\frac{Q_4}{Q_{\alpha+1}N^{\frac{1}{2}}}<\frac{1}{Q^2_{\alpha+1}}\le\frac{1}{\q},
\end{equation}
откуда $\xi=0$.

Наоборот, пусть теперь $Q_{\alpha+1}\ge\frac{1}{Q_4}N^{\frac{1}{2}}$ (\textbf{случай 2.2}). Тогда, ввиду (\ref{27norm}) и (\ref{28norm}), при $l=l^{(1)},l^{(2)}$ выполнено неравенство
\begin{equation} 
\label{87norm}
|l|\le Q_{\beta+1}\le\frac{3Q_3N^{\frac{1}{2}}}{Q_{\alpha+1}}\le3Q_7.
\end{equation}
Поэтому, в обозначениях (\ref{43norpm}), имеет место оценка 
\begin{equation} 
\label{88norm}
|l^{(1)}-l^{(2)}|\le6Q_7<7Q_7=T_1.
\end{equation}
Однако сравнение в (\ref{43norm}) показывает, что $l^{(1)}\equiv l^{(2)}\pmod {T_1}$. Вместе с неравенством (\ref{88norm}) это означает, что $l^{(1)}= l^{(2)}$. Возвращаясь теперь к неравенству (\ref{85norm}), заменим в нем $F(x_2)$ на первую строку минимума в (\ref{52norm}), равную, ввиду (\ref{81norm}), числу $\frac{1}{N}150A^3x_2$. Тогда получим:
$$\left\|\xi \right\|\le\frac{150A^3}{N}<\frac{1}{q^{(1)}q^{(2)}}\le \frac{1}{\q}$$
ввиду неравенства в (\ref{17norm}), откуда снова получаем равенство $\xi=0$. Таким образом, $\frac{a^{(1)}}{q^{(1)}}=\frac{a^{(2)}}{q^{(2)}}$. Отсюда и из (\ref{81norm}) следует равенство (\ref{65norm}). Лемма доказана.

\begin{Th}\label{t8.1}
 Пусть число $M^{(1)}$ определено любым из равенств (\ref{68norm}) или (\ref{71norm}). Тогда в условиях теоремы 6.1 имеет место оценка  
\begin{equation} 
\label{90norm}
\sigma_{N,Z}\ll \left( M^{(1)}\right)^{1+\epsilon_0}\left   |\Omega^{(1)}\right|^{\frac{1}{2}}
\sum_{g_3\in \Omega^{(3)}}\left|\M(g_3)\right|^{\frac{1}{2}}.
\end{equation}
\end{Th}

Доказательство. Из определения числа  $M^{(1)}$ следует, что неравенство (\ref{54norm}) выполнено по крайней мере для всех достаточно больших значений $N$. Следовательно, условия леммы \ref{8.3} выполнены. Значит, имеет место неравенство (\ref{55norm}). Отсюда и из включения (\ref{69norm}), доказанного в леммах \ref{8.4} и \ref{8.6}, следует неравенство (\ref{90norm}) (если учесть определение числа $H$ в  (\ref{norkkm})). Теорема доказана.

Тем самым задача об оценке тригонометрической суммы сведена к определению мощности множества $\M(g_3)$.

\section{Исследование и применение множества $\M(g_3)$ }

В этом параграфе будет получено значение мощности множества $\M(g_3)$, что позволит вывести полезные следствия из соотношения (\ref{90norm}). Всюду далее число $M^{(1)}$ определено одним из двух вариантов --- как в теореме \ref{t8.1}.

Напомним обозначения (\ref{lorm}) и (\ref{24norm}) и при $\beta\ge0$ положим:
\begin{equation} 
\label{91norm}
M^{(2)}_{\beta}=\frac{1}{Q_3}\left(M^{(1)}\right)^{-2\epsilon_0}\left(Q_{\beta}\right)^{\frac{1}{2}-2\epsilon_0},
\end{equation}
\begin{equation} 
\label{92norm}
M^{(2)}=
\left\{
              \begin{array}{ll}
                M^{(2)}_{\beta},{}\ \mbox{если  }M^{(2)}_{\beta}\ge Q_1,\\
               1, \mbox{--- в противном случае.}
              \end{array}
\right.
\end{equation}

В следующих ниже леммах \ref{9.1} --- \ref{9.3} свойства (\ref{normyy}) --- (\ref{8norm}) продолжаем считать выполненными.

\begin{Le}\label{9.1}  Пусть число $M^{(2)}$ определено соотношением (\ref{92norm}). Тогда при любом $g_3~\in~\Omega^{(3)}$ для каждого элемента множества $\M(g_3)$ выполняется равенство
\begin{equation} 
\label{93norm}
g^{(1)}_2=g^{(2)}_2.
\end{equation}\end{Le}

Доказательство.  Если $M^{(2)}=1$, то $\Omega^{(2)}=\{E\}$, так что утверждение леммы было бы очевидно. В частности, если $Q_{\beta}=0$, то снова $M^{(2)}=1$, и лемма была бы доказана. 

Пусть теперь утверждение леммы не имеет места. Тогда 
\begin{equation} 
\label{94nnorm}
M^{(2)}=M^{(2)}_{\beta}\ge Q_1>1,{}\ {}\ |l^{(1)}|,|l^{(2)}|\ge Q_{\beta}>0,
\end{equation}
и найдется четверка 
\begin{equation} 
\label{94norm}
  \left(\widetilde{g}^{(1)},\widetilde{g}^{(2)},\Theta^{(1)},\Theta^{(2)}\right) \in  \M(g_3),
\end{equation}
	для которой $g^{(1)}_2\neq g^{(2)}_2$. Последнее неравенство, ввиду взаимной однозначности отображения (\ref{200norm}), имеет особую интерпретацию в обозначениях (\ref{11norm}). Именно, не только выполнено соотношение $\frac{x_1}{x_2}\neq\frac{y_1}{y_2}$. Более того, если последовательность $D$ (возможно, пустая) --- максимальная общая часть последовательностей неполных частных для цепных дробей 
\begin{equation} 
\label{95norm}
\frac{x_1}{x_2}=[D,T],{}\ {}\ \frac{y_1}{y_2}=[D,W],
\end{equation}
то последовательность $D$ не длиннее, чем те последовательности, которым соответствуют матрицы $g^{(1)}_2$ и $g^{(2)}_2$ в смысле отображения (\ref{mattrix-deff}). Следовательно,  из (\ref{6norm}) имеем оценку:
\begin{equation} 
\label{96norm}
\langle D\rangle\le\max\limits_{g_2\in\Omega^{(2)}}\|g_2\|\le
73A^2   M^{(2)}{ \left  (M^{(1)}M^{(2)}\right)}^{2\epsilon_0}.
\end{equation}
Кроме того, длины цепных дробей (\ref{95norm}) совпадают по построению ансамбля. Следовательно, выполнены условия леммы \ref{5.1}, применяя которую, получаем:
\begin{equation} 
\label{97norm}
\left|\frac{x_1}{x_2}-\frac{y_1}{y_2}\right|\ge \frac{1}{(2A)^{4}\langle D\rangle^2}
\ge \frac{1}{(2A)^{4}\left(73A^2   M^{(2)}{ \left  (M^{(1)}M^{(2)}\right)}^{2\epsilon_0} \right)^2}\ge 
\frac{(5A)^{-8}\left(M^{(1)}\right)^{-4\epsilon_0}}
{\left(M^{(2)}\right)^{2+4\epsilon_0}}
\end{equation}
ввиду (\ref{96norm}). Подставляя сюда выражение для $M^{(2)}=M^{(2)}_{\beta}$ из (\ref{91norm}) и применяя оценку $ Q_3\ge {\left(7A \right)}^7,$ получаем:
\begin{equation} 
\label{98norm}
\left|\frac{x_1}{x_2}-\frac{y_1}{y_2}\right|\ge \frac{{\left(7A \right)}^{7{(2+4\epsilon_0)}}}{(5A)^{8}Q_{\beta}}>\frac{2(9A)^{5}}{Q_{\beta}}.
\end{equation}
С другой стороны, применяя при $i=1,2$ неравенство (\ref{66norm}), получаем:
$$\left|\frac{x_1}{x_2}-\frac{y_1}{y_2}\right|\le \frac
{y_1\left| x_2l^{(1)}-y_2 l^{(2)} \right|+y_2\left| x_1l^{(1)}-y_1 l^{(2)} \right|}{x_2y_2l^{(1)}}\le\frac{(9A)^{5}}{l^{(1)}}\left( \frac{y_1}{y_2}+
\frac{x_1}{x_2}\right)\le\frac{2(9A)^{5}}{Q_{\beta}}$$
(напомним, что $|l|^{(1)}>0$ ввиду (\ref{94nnorm})), что противоречит неравенству (\ref{98norm}). К этому противоречию привело предположение о том, что $g^{(1)}_2\neq g^{(2)}_2$. Следовательно, $g^{(1)}_2= g^{(2)}_2$. Лемма доказана.

Для $\alpha\ge0$ положим: 
\begin{equation} 
\label{100norm}
M^{(4)}_{\alpha}=\frac{1}{Q_3}\left( Q_{\alpha}\right)^{\frac{1}{2}},
\end{equation}
\begin{equation} 
\label{101norm}
M^{(4)}=
\left\{
              \begin{array}{ll}
                M^{(4)}_{\alpha},{}\ \mbox{если  }M^{(4)}_{\alpha}\ge Q_1,\\
               1, \mbox{--- в противном случае.}
              \end{array}
\right.
\end{equation}
Всюду далее числа $M^{(2)}$ и $M^{(4)}$ определены соотношениями (\ref{92norm}) и (\ref{101norm}), соответственно. 

\begin{Le}\label{9.2}  При любом $g_3\in\Omega^{(3)}$ для каждого элемента множества $\M(g_3)$ выполняются равенства (\ref{93norm}) и
\begin{equation} 
\label{102normlttt}
\widetilde{g}^{(1)}_4=\widetilde{g}^{(2)}_4.
\end{equation}\end{Le}

  Доказательство. Если $M^{(4)}=1$, то $\Omega^{(4)}=\{E\}$, так что равенство (\ref{102normlttt}) было бы очевидно. В частности,  если $Q_{\alpha}=1$, то снова $M^{(4)}=1$, и лемма была бы доказана.
  
  Пусть теперь утверждение леммы не имеет места. Тогда $M^{(4)}=M^{(4)}_{\alpha}\ge Q_1~>~1, \\  \q~ \ge ~ Q_{\alpha}> 1$ и, следовательно, $a^{(1)}>0$.
  
  Равенство (\ref{93norm}), согласно которому $g^{(1)}_2=g^{(2)}_2=g_2$, уже доказано в предыдущей лемме. Следовательно, сравнения в (\ref{67norm}) могут быть записаны в виде
$$\left(g_2g_3a^{(1)}\widetilde{g}^{(1)}_4\equiv g_2g_3a^{(2)}\widetilde{g}^{(2)}_4 \right)_{1,2}\pmod\q,$$
где индексы ``1,2'' внизу означают, как всегда, выполнение сравнения по обеим координатам. Отсюда, в виду равенства $\det g_2g_3=1$, получаем: 
\begin{equation} 
\label{103norm}
\left(a^{(1)}\widetilde{g}^{(1)}_4\equiv a^{(2)}\widetilde{g}^{(2)}_4 \right)_{1,2}\pmod\q.
\end{equation}
Положим $\widetilde{g}^{(1)}_4=(u_1,u_2)^t,{}\ {}\ \widetilde{g}^{(2)}_4=(v_1,v_2)^t,$ тогда сравнения (\ref{103norm}) перепишутся в виде
$$a^{(1)}u_1\equiv a^{(2)}v_1 \pmod\q,$$
$$a^{(1)}u_2\equiv a^{(2)}v_2 \pmod\q.$$
Отсюда следует цепочка сравнений:
\begin{equation} 
\label{102norm}
\left( a^{(1)}u_2\right)v_1\equiv  a^{(2)}v_2v_1= 
\left( a^{(2)}v_1\right)v_2\equiv \left( a^{(1)}u_1\right)v_2\pmod\q.
\end{equation}
Но числа $a^{(1)}$ и $\q$ взаимно просты как числитель и знаменатель цепной дроби (и при этом $a^{(1)}>0$, как уже было сказано выше). Следовательно, сокращая начало и конец цепочки (\ref{102norm}) на $a^{(1)}$, получаем:
\begin{equation} 
\label{105norm}
u_2v_1-u_1v_2\equiv0\pmod\q.
\end{equation}
С другой стороны, 
\begin{equation} 
\label{106norm}
|u_2v_1-u_1v_2|\le u_2v_2 \le \left(\max\limits_{g\in\Omega^{(4)}}\|g\|\right)^2.
\end{equation}
Применяя к окончанию соотношения (\ref{106norm}) верхнюю оценку из
 неравенства (\ref{7norm}), получаем: 
 $$|u_2v_1-u_1v_2|\le\left(73A^2M^{(4)}\right)^2\le
 (9A)^4\left(M^{(4)}\right)^2. $$
 Подставляя в последнее неравенство значение $M^{(4)}=M^{(4)}_{\alpha}$ из (\ref{100norm}), получаем:
\begin{equation} 
\label{107norm}
|u_2v_1-u_1v_2|<Q_{\alpha}\le\q.
\end{equation}
ввиду нижней оценки из неравенства \ref{26norm}. Соотношения (\ref{105norm}) и (\ref{107norm}) показывают, что 
$$u_2v_1=u_1v_2.$$
Отсюда, вследствие несократимости дробей $\frac{u_1}{u_2},{}\ \frac{v_1}{v_2}$, получаем: $u_1=v_1,{} \ u_2=v_2,$ или $\widetilde{g}^{(1)}_4=\widetilde{g}^{(2)}_4,$ так что равенство (\ref{102normlttt}) доказано. Лемма доказана.

\begin{Le}\label{9.3}  При любом $g_3\in\Omega^{(3)}$ имеет место равенство
\begin{equation}
\label{108norm}
\left |\M(g_3)\right|=\left   |\Omega^{(2)}\right|\left   |\Omega^{(4)}\right||Z|.
\end{equation}\end{Le}

Доказательство. Пересчитаем четверки элементов (\ref{94norm}), пользуясь  предыдущей леммой. Для этого выберем и фиксируем матрицы ${g}^{(1)}_2={g}^{(2)}_2$ и ${g}^{(1)}_4={g}^{(2)}_4$, определяющие элементы $\widetilde{g}^{(1)}=\widetilde{g}^{(2)}$, одним из $\left   |\Omega^{(2)}\right|\left   |\Omega^{(4)}\right|$ способов --- это первые два множителя в (\ref{108norm}). Выберем также число $\Theta^{(2)} \in Z,$ как в  (\ref{50norm}), $|Z|$ способами --- это третий и последний множитель в (\ref{108norm}). Заметим теперь, что из равенства $\widetilde{g}^{(1)}=\widetilde{g}^{(2)}$ (вытекающего из предыдущей леммы) следует, что $x_1=y_1.$ Поэтому, согласно замечанию \ref{z8.1}, выполняются равенства $a^{(1)}=a^{(2)}$ и $l^{(1)}=l^{(2)}$. Подстановка последних в (\ref{50norm}) приводит к соотношению $\Theta^{(1)}=\Theta^{(2)}$, завершающему вывод формулы (\ref{108norm}). Лемма доказана.

Далее выполнение условий теоремы \ref{t6.1} заранее не предполагается. Напомним обозначения (\ref{29norm}) и (\ref{30norm}).

\begin{Th}\label{t9.1}
 Для всякого числа $M^{(1)}$,  определенного любым из равенств (\ref{68norm}) или (\ref{71norm}), из неравенства 
\begin{equation} 
\label{109norm}
M^{(1)}M^{(2)}_{\beta}M^{(4)}_{\alpha}\le N
\end{equation}
следует оценка
\begin{equation} 
\label{110norm}
\sigma_{N,Z}\ll
   \frac{\left   |\Omega_N\right||Z|}
  {(\overline{Q}_{\alpha}\overline{Q}_{\beta}|Z|)^{\frac{1}{2}}}
   \left(  {M^{(1)}Q^{\frac{1}{2}}_{\alpha+1}Q^{\frac{1}{2}}_{\beta+1}}\right)^{ \gamma+8\epsilon_0}.
\end{equation}
\end{Th}

Доказательство. Ввиду неравенства (\ref{109norm}), выполнено условие (\ref{yorm}) теоремы \ref{t6.1}, согласно которой найдется разложение ансамбля (\ref{1norm}) со свойствами (\ref{2norm}) --- (\ref{8norm}). Следовательно, для оценки величины $\sigma_{N,Z}$ имеем право применять как теорему \ref{t8.1}, так и лемму \ref{9.3}, утверждения которых приводят к неравенству:
$$
\sigma_{N,Z}\ll 
\left(M^{(1)}\right)^{1+2\epsilon_0}\left   |\Omega^{(1)}\right|^{\frac{1}{2}}
\sum_{g_3\in \Omega^{(3)}}
\left   |\Omega^{(2)}\right|^{\frac{1}{2}}
\left   |\Omega^{(4)}\right|^{\frac{1}{2}}
|Z|^{\frac{1}{2}}.$$
Отсюда, применяя соотношения (\ref{1norm}) --- (\ref{2norm}), получаем:
\begin{equation} 
\label{111norm}
\sigma_{N,Z}\ll 
\frac{\left(M^{(1)}\right)^{1+2\epsilon_0}\left   |\Omega_N\right|
|Z|^{\frac{1}{2}}}
{\left   |\Omega^{(1)}\Omega^{(2)}\right|^{\frac{1}{2}}
\left   |\Omega^{(4)}\right|^{\frac{1}{2}}
}\ll
\frac{\left(M^{(1)}\right)^{1+2\epsilon_0}\left   |\Omega_N\right|
|Z|^{\frac{1}{2}}}
{\left( {M^{(1)}M^{(2)}M^{(4)}}\right)^{\delta_{\A}-\epsilon_0}}.
\end{equation}
Подставляя в неравенство (\ref{111norm}) определения величин $M^{(2)}$  и $M^{(4)}$, получаем оценку (\ref{110norm}). Теорема доказана.

\section{Доказательство теорем 2.1 и 2.2}
Согласно теореме \ref{t7.1}, для доказательства обеих теорем из заголовка параграфа достаточно доказать неравенство (\ref{45norm}) при $c_2=5$. Будем выводить оценку (\ref{45norm}) из неравенства (\ref{110norm}) теоремы \ref{t9.1}. 

 Для начала определим  число $M^{(1)}$ равенством (\ref{71norm}), тогда проверка неравенства (\ref{109norm}) легко получается применением замечания \ref{z7.2}. Следовательно, выполнено условие теоремы \ref{t9.1}, применение которой приводит к оценке
\begin{equation} 
\label{112norm}
\sigma_{N,Z}\ll
   \frac{\left   |\Omega_N\right||Z|}
  {(\overline{Q}_{\alpha}\overline{Q}_{\beta}|Z|)^{\frac{1}{2}}}
   \left(  N^{\frac{1}{2}}Q_{\alpha+1}Q_{\beta+1}\right)
   ^{ \gamma}N^{12\epsilon_0}.
\end{equation}
Пусть (\textbf{случай 1}) выполнено неравенство
\begin{equation} 
\label{113norm}
N<(Q_{\alpha+1})^{\frac{5}{2}}(Q_{\beta+1})^{\frac{3}{2}}.
\end{equation} 
Тогда, ввиду (\ref{113norm}), имеет место оценка
 $$N^{\frac{1}{2}}Q_{\alpha+1}Q_{\beta+1}
 <(Q_{\alpha+1})^{\frac{9}{4}}(Q_{\beta+1})^{\frac{7}{4}}<
 (Q_{\alpha+1})^{\frac{5}{2}}(Q_{\beta+1})^{\frac{5}{2}},$$
 подстановка которой в неравенство (\ref{112norm}) приводит к нужной оценке  (\ref{45norm}) с $c_2=5$. 
 
 Если же неравенство (\ref{113norm}) не выполнено (\textbf{случай 2}), то, следовательно, имеет место противоположное неравенство 
 $$(Q_{\alpha+1})^{\frac{5}{2}}(Q_{\beta+1})^{\frac{3}{2}}\le N.$$
 В этом случае для числа $M^{(1)}$, определенного равенством (\ref{68norm}), выполнено неравенство (\ref{109norm}). Следовательно, имеем право применить теорему \ref{t9.1}, которая при подстановке такого значения $M^{(1)}$ приводит  к оценке (\ref{45norm}) с $c_2=5$. Доказательство теорем закончено.

\end{document}